# A finite coefficient theory for shape identification in tilted generalized gamma families

Ran Huang

Centre of Biomedical Engineering and Technology, Yiwu Research Institute of Fudan University, Yiwu 322000, China

Correspondence: huangran@fudan.edu.cn

## Abstract

A continuous shape parameter enters both a radial density and the normalization of its transform coefficients. This coupling allows changes in shape to be compensated by changes in scale and exponential tilt. We study the coefficients $b_k = e^a \mathbb{E}[Y^k]/[k!\,(d)_k]$ for densities proportional to $y^{d-1}\exp(-\alpha y^q + \beta y)$, with known integer $q \geq 2$, $d, \alpha > 0$, and unrestricted real $a, \beta$. For $q = 2$, the first four coefficients determine all parameters globally, whereas three coefficients leave an analytic curve of solutions. We characterize the exact image of the two shape ratios and prove that two alternating corners decide whether an entire coefficient box is realizable and give its exact shape range. For every fixed $q$, the first $q(q+1)$ coefficients give global identification and an everywhere full-rank map. A finite-difference criterion yields generic identification from $2q + 1$ coefficients, including every zero-tilt model against the full tilted family. For $q = 3$, an exact rational certificate constructs two distinct shapes with identical first four coefficients and separated fifth coefficients. The results distinguish exact identification, realizability, and sensitivity for a parameter-dependent coefficient normalization.

**Keywords:** structural identifiability; generalized gamma family; normalized moments; parabolic cylinder functions; validated computation

## 1 Introduction

A radial response can change because its characteristic scale changes, because the potential redistributes mass, or because the radial measure itself changes. A small set of response coefficients may represent these changes in the same way. The resulting ambiguity concerns the observation map as well as the underlying distribution: the dimension or shape parameter may occur in the kernel that converts radial moments into observable coefficients. Determining a probability density from its ordinary moments and determining its shape from such coefficients are therefore different inverse problems.

Consider a positive radial variable with density proportional to $y^{d-1}\exp(-\alpha y^q + \beta y)$. The parameter $d$ controls the weight of the measure near the origin, $\alpha$ sets a scale, and $\beta$ tilts the distribution. When $\beta = 0$, this is a generalized gamma density. For quadratic and cubic potentials in $y$, the corresponding even weights after $y = x^2$ belong to the generalized Freud families. Their moments, Pearson equations, and connections with special functions and orthogonal polynomials have an established theory [1,2]. We use these existing models to ask a different question: how many consecutive coefficients of a transform with the unknown normalization $(d)_k$ distinguish all four parameters, including an amplitude?

For a minimal exponential family, strict convexity of the log-partition function identifies the natural parameters from the means of its sufficient statistics [3]. In the present density those statistics can be taken as $\log Y$, $Y$, and $Y^q$. Our data instead divide the ordinary power moments by $(d)_k$. Equality of the observed coefficients for two shapes therefore permits unequal ordinary moments. The convexity argument for the usual mean map does not by itself settle the present inverse problem.

Two further bodies of work provide relevant tools. For the quadratic potential, the normalized moments are expressible through parabolic cylinder functions. Sharp bounds and monotonicity properties of their ratios and double ratios are known [4,5], with a recent account in [6]. These results control the tilt variable at fixed order. Joint recovery of the order and the argument requires an additional comparison along a curve on which one double ratio is fixed. On the algebraic side, reconstructing functions from moments by converting differential equations into recurrences is an established method for piecewise D-finite and holonomic densities [7,8]. Holonomic methods also support estimation in exponential-polynomial families [9]. Our recurrence argument operates on coefficients whose denominator changes with the unknown shape, and gives explicit global bounds for this observation map.

The quadratic case admits a complete answer. After removing amplitude and scale, two shape ratios form global analytic coordinates on an explicitly characterized open region. The main step compares models at nearby shapes through a beta multiplication that exactly compensates the Pochhammer factors. The infinitesimal change has a strictly negative third derivative after its polynomial part is removed. Orthogonality then fixes the sign of the third response while the first two responses remain constant. The same image geometry yields a sharp statement for independent coefficient intervals: two specific corners decide whether every datum in the box is realizable, and they give the exact minimum and maximum shape.

For a general fixed integer potential degree, a finite comparison of two Pearson recurrences gives global identification from $q(q+1)$ coefficients, with full differential rank everywhere. A second argument yields identification from $2q+1$ coefficients whenever a data-dependent finite difference is nonzero. We prove that this condition holds throughout the zero-tilt slice and on an open dense subset of the full parameter domain. A shorter elimination limits the number of preimages to two using $q+3$ coefficients when $q \geq 3$. Finally, a certified cubic example shows that the quadratic four-coefficient theorem does not extend unchanged to higher potentials.

The contribution is the global geometry and finite identification theory of this specific coefficient map. The density family, beta identities, moment positivity, and recurrence method are classical ingredients. We give the proofs needed to separate their roles. The general coefficient bound is sufficient rather than optimal, and the cubic example separates a particular pair with its fifth coefficient. Neither statement supplies a global five-coefficient theorem for the cubic family.

## 2 Model and observation map

Fix an integer $q \geq 2$. Write

$$Z_q(d,\alpha,\beta) = \int_0^\infty y^{d-1}\, e^{-\alpha y^q + \beta y}\, dy, \qquad p(y) = \frac{y^{d-1} e^{-\alpha y^q + \beta y}}{Z_q(d,\alpha,\beta)}, \tag{1}$$

where $d, \alpha > 0$ and $\beta \in \mathbb{R}$. The parameter domain is $\Theta = \mathbb{R} \times (0,\infty)^2 \times \mathbb{R}$, with coordinates $(a, d, \alpha, \beta)$. Our observations are

$$b_k = \frac{e^a m_k}{k!\,(d)_k}, \qquad m_k = \mathbb{E}[Y^k], \qquad (d)_k = \prod_{j=0}^{k-1}(d+j), \quad (d)_0 = 1. \tag{2}$$

The phrase "first $K$ coefficients" means $b_0, \dots, b_{K-1}$ and includes $b_0 = e^a$. If amplitude is externally known, this coefficient carries no additional shape information; our counts retain it to keep the full observation map explicit.

These data are the Taylor coefficients of

$$\mathcal{T}(z) = e^{a}\mathbb{E}[{}_0F_1(;d;zY)] = \sum_{k=0}^{\infty} b_k\, z^k. \qquad (3)$$

The expectation and series may be interchanged uniformly on compact $z$ sets. Indeed, with $c = \min(d, 1)$, $(d)_k \geq c^k k!$, so the absolute kernel series is bounded by $\exp\left(2\sqrt{|z|Y/c}\right)$, which is integrable under (1). Thus (3) is entire. For integer spatial dimension $D = 2d$, the normalized spherical Fourier kernel at radius $R = 2\sqrt{Y}$ is ${}_0F_1(;d;-t^2Y)$; the sign-adjusted even Fourier coefficients are then (2). The familiar normalized Bessel kernel and its beta dimension relation are classical [10]. For general $d > 0$, (3) defines the analytic observation model directly. No geometric realization of every noninteger dimension is assumed.

The specification of (3) matters experimentally. One may regard the $b_k$ as measured response coefficients without knowing $d$ in advance. If only samples of $Y$ are available, however, dividing their empirical moments by an unknown $(d)_k$ does not create parameter-independent observations. A statistical estimation problem from such samples would require its own observation model.

Set

$$c_k = \frac{b_k b_0^{k-1}}{b_1^k}, \qquad c_0 = c_1 = 1, \qquad u = c_2, \quad v = c_3. \qquad (4)$$

This normalization replaces $Y$ by $Y' = b_0 Y/b_1$ and removes amplitude. The transformed parameters are $\alpha' = \alpha(b_1/b_0)^q$ and $\beta' = \beta b_1/b_0$. It is reversible once $b_0, b_1$ are retained. In particular, it preserves the number of full parameter preimages of any fixed datum.

All coefficient maps in this paper are real analytic. In a sufficiently small complex neighborhood of a real parameter point, $\mathrm{Re}d$ and $\mathrm{Re}\alpha$ have positive lower bounds. At zero, powers of $y$ times powers of $|\log y|$ are integrable; at infinity, $e^{-cy^q+Cy}$ dominates each required derivative. The normalization integral is positive at the real point and remains nonzero after shrinking the neighborhood. This also justifies all finite parameter differentiations used below.

The principal counts are collected in Table 1. "Generic" means an open dense subset of $\Theta$ whose data each have a unique preimage in the entire $\Theta$.

| Potential | Coefficients | Conclusion |
|---|---|---|
| q = 2 | 3 | An analytic curve of preimages for each realizable datum |
| q = 2 | 4 | Global uniqueness and exact realizability region |
| Any fixed q ≥ 2 | q(q + 1) | Global uniqueness and rank four everywhere |
| Any fixed q ≥ 2 | 2q + 1 | Generic global uniqueness; all zero-tilt models included |
| Any fixed q ≥ 3 | q + 3 | At most two full parameter preimages |
| q = 3 | 4 | A certified pair with different shapes and identical data |

Table 1. Coefficient counts include amplitude $b_0$. The general sufficient bounds are not asserted to be minimal. The quadratic four-coefficient theorem is stronger than the general bounds at $q = 2$.

## 3 Global identification for the quadratic potential

For $q = 2$, scale out $\alpha$ and define

$$J_k(d,\lambda) = \int_0^\infty t^{d+k-1}\, e^{-t^2+\lambda t}\, dt, \qquad \lambda = \beta/\sqrt{\alpha}, \qquad r = J_1/J_0. \qquad (5)$$

The two scale-invariant ratios are

$$u(d,\lambda) = \frac{dJ_0J_2}{2(d+1)J_1^2}, \qquad v(d,\lambda) = \frac{d^2J_0^2J_3}{6(d+1)(d+2)J_1^3}. \qquad (6)$$

Introduce

$$A(\lambda) = \int_0^\infty e^{-t^2+\lambda t}\, dt, \qquad U_0(\lambda) = \frac{A'(\lambda)}{2A(\lambda)^2}, \qquad W_0(\lambda) = \frac{A''(\lambda)}{12A(\lambda)^3}. \qquad (7)$$

We prove below that $U_0$ decreases analytically from $1/2$ to $0$. For $0 < u < 1/2$, let $\lambda_0(u) = U_0^{-1}(u)$ and put

$$v_-(u) = \frac{u^2}{3(1-u)}, \qquad v_+(u) = W_0\big(\lambda_0(u)\big), \qquad (8)$$
$$\Omega = \{(u,v)\colon 0 < u < 1/2,\ v_-(u) < v < v_+(u)\}.$$

**Theorem 1.** The map $(d,\lambda) \mapsto (u,v)$ is a real analytic diffeomorphism from $(0,\infty)\times\mathbb{R}$ onto $\Omega$, and its Jacobian in the column order $(d,\lambda)$ is strictly negative. A positive four-vector $(b_0, b_1, b_2, b_3)$ is realizable if and only if the ratios in (4) lie in $\Omega$. Its unique full parameters are

$$a = \log b_0, \qquad (d,\lambda) = H^{-1}(u,v), \qquad \alpha = \left(\frac{b_0 r}{db_1}\right)^2, \qquad \beta = \lambda\sqrt{\alpha}. \qquad (9)$$

For any fixed positive $(b_0, b_1, b_2)$ with $0 < u < 1/2$, the preimages form one analytic curve parametrized by $0 < d < L(u)$, where $L(u) = 2u/(1-2u)$. Along this curve $b_3$ decreases strictly with $d$. Consequently, four is the minimal number of consecutive coefficients for global identification in this family.

### 3.1 Fixed shape monotonicity from a moment determinant

Integration by parts gives

$$2J_{k+2} = \lambda J_{k+1} + (d+k)J_k. \qquad (10)$$

Let $R_j = J_{j+1}/J_j$. Subtracting $\lambda/2$ times the second column and $(d+j)/2$ times the first column from the third column of the $3\times 3$ moment matrix yields

$$\det\big[J_{j+i+k}\big]_{i,k=0}^2 = \frac{J_jJ_{j+1}J_{j+2}}{2}\big(2R_{j+1} - R_j - R_{j+2}\big) > 0. \qquad (11)$$

Strict positivity follows because this matrix is the Gram matrix of $1, t, t^2$ under a positive density on the whole half-line. Therefore the sequence $R_j$ is strictly discretely concave. Since $\partial_\lambda J_j = J_{j+1}$,

$$\partial_\lambda \log u = R_0 + R_2 - 2R_1 < 0. \qquad (12)$$

This fixed-order monotonicity is also a consequence of the established parabolic cylinder double-ratio theory [4,5]. Equation (11) gives the short self-contained argument needed here.

For fixed $d$, scaling $t = x/|\lambda|$ as $\lambda \to -\infty$ gives gamma moments and $u \to 1/2$. For $\lambda \to +\infty$, completing the square at $t = \lambda/2$ gives relative concentration and $u \to d/[2(d+1)]$. The convergence in the latter limit is uniform for $d$ in a compact subset of $(0,\infty)$: split $|t - \lambda/2| \le \lambda/4$ from its complement, expand the bounded relative power on the central interval, and bound the complement with the Gaussian tail; the near-zero part is integrable uniformly because $d$ stays positive. Thus, at fixed $u \in (0,1/2)$, there is exactly one analytic solution $\lambda_u(d)$ for $0 < d < L(u)$.

### 3.2 Comparison between shapes

Let $T_d$ have the normalized density from (5) along this fixed-$u$ curve, and set $X_d = dT_d/r(d)$. Its density is proportional to $x^{d-1}e^{-\alpha_d x^2+\beta_d x}$, where $\alpha_d = (r/d)^2$ and $\beta_d = \lambda_u(d)r/d$. Its first three moments obey

$$\mathbb{E}X_d = d, \quad \mathbb{E}X_d^2 = 2u\,d(d+1), \quad \mathbb{E}X_d^3 = 6V(d)(d)_3, \qquad V(d) = v\big(d, \lambda_u(d)\big). \qquad (13)$$

Fix $d > 0$. For small $\varepsilon > 0$, take $B_\varepsilon \sim \mathrm{Beta}(d, \varepsilon)$ independently of $X_{d+\varepsilon}$ and define $Y_\varepsilon = B_\varepsilon X_{d+\varepsilon}$. The classical beta moment identity gives

$$\mathbb{E}Y_\varepsilon^k = \frac{(d)_k}{(d+\varepsilon)_k}\mathbb{E}X_{d+\varepsilon}^k. \qquad (14)$$

The moments of orders zero, one, and two are constant in $\varepsilon$, while the third equals $6(d)_3V(d+\varepsilon)$. This use of the beta kernel is consistent with the classical Sonine dimension relation [10]. Its role here is to compare the changing normalization in (2).

Write $\alpha_\varepsilon = \alpha_{d+\varepsilon}$ and $\beta_\varepsilon = \beta_{d+\varepsilon}$. A direct product-density calculation gives the density of $Y_\varepsilon$ as

$$q_\varepsilon(y) = \frac{\Gamma(d+\varepsilon)y^{d-1}e^{-\alpha_\varepsilon y^2+\beta_\varepsilon y}}{\Gamma(d)Z_2(d+\varepsilon,\alpha_\varepsilon,\beta_\varepsilon)}F_\varepsilon(y), \qquad (15)$$

$$F_\varepsilon(y) = \frac{1}{\Gamma(\varepsilon)}\int_0^\infty s^{\varepsilon-1}\, e^{-\alpha_\varepsilon s^2+(\beta_\varepsilon-2\alpha_\varepsilon y)s}\, ds. \qquad (16)$$

Subtract the integral representation of $\Gamma(\varepsilon)$ inside (16). Since $1/\Gamma(\varepsilon) = \varepsilon + O(\varepsilon^2)$, the first-order nonpolynomial term is

$$G(y) = \int_0^\infty \frac{e^{-\alpha s^2+(\beta-2\alpha y)s} - e^{-s}}{s}\, ds, \qquad \alpha = \alpha_d, \quad \beta = \beta_d. \qquad (17)$$

For completeness, the limiting operation is valid in each polynomially weighted $L^1$ space used here. On a small parameter interval, choose $0 < a \le \alpha_\varepsilon \le A$ and $|\beta_\varepsilon| \le b$. On $0 < s < 1$, the difference of exponentials in the numerator is bounded in absolute value by $C(1+y)s$; one uses $|e^{-2\alpha_\varepsilon ys} - 1| \le 2Ays$. On $s \ge 1$, a bound is $e^{-as^2+bs} + e^{-s}$. Thus the subtracted integral is at most $C(1+y)$ uniformly for small $\varepsilon$. The remaining density factor supplies a common integrable bound of the form $Cy^{d-1}(1+y)^M e^{-ay^2+by}$. The difference quotients in (15) are consequently integrable after multiplication by $y^k$, $0 \le k \le 3$.

The derivative at zero therefore has the form $p_d(y)\psi(y)$, with

$$\psi(y) = C_0 - \alpha_d{}'y^2 + \beta_d{}'y + G(y), \qquad \int y^k\psi(y)p_d(y)\,dy = 0 \quad (k = 0,1,2). \qquad (18)$$

Differentiation of (17), justified by its Gaussian tail, gives

$$\psi'''(y) = -8\alpha^3\int_0^\infty s^2\, e^{-\alpha s^2+(\beta-2\alpha y)s}\, ds < 0. \qquad (19)$$

Let $P_3$ be the monic cubic orthogonal polynomial for $p_d(y)\,dy$. Its three zeros are distinct and positive. Indeed, if there were fewer than three sign changes in the support, multiplying its sign-change factors would give a polynomial of degree at most two whose product with $P_3$ has a fixed nonzero sign, contradicting orthogonality. Interpolate $\psi$ at these zeros by a quadratic $Q_2$. The Lagrange remainder and (19) give $(\psi - Q_2)P_3 < 0$ away from the zeros. Using (18) and orthogonality,

$$6(d)_3 V'(d) = \int y^3 \psi p_d = \int P_3(\psi - Q_2) p_d < 0. \qquad (20)$$

The curve is analytic, so the right derivative obtained from beta multiplication is its ordinary derivative. We have proved strict decrease along every fixed-$u$ curve. In particular,

$$\det DH = u_d v_\lambda - u_\lambda v_d = -u_\lambda V'(d) < 0. \qquad (21)$$

### 3.3 The two boundaries and completion of the proof

First consider $d \downarrow 0$. In a neighborhood of any finite real $\lambda$,

$$J_0(d,\lambda) = \frac{1}{d} + \mathcal{B}(d,\lambda), \qquad \mathcal{B} = \int_0^1 t^{d-1}\left(e^{-t^2+\lambda t} - 1\right) dt + \int_1^\infty t^{d-1}\, e^{-t^2+\lambda t}\, dt. \qquad (22)$$

Here $\mathcal{B}$ is jointly analytic for $d$ near zero. The expressions (6) extend analytically to $d = 0$ and become $U_0, W_0$ in (7). To see that $U_0{}' < 0$, symmetrize $AA'' - 2(A')^2$ and substitute $s + t = x$, $t/x = p$. This gives

$$AA'' - 2(A')^2 = \int_0^\infty x^3\, e^{-x^2+\lambda x} \int_0^1 \left(\frac{1}{2} - 3p(1-p)\right) e^{2x^2 p(1-p)}\, dp\, dx < 0. \qquad (23)$$

The inner integral is $-3\mathrm{Cov}\left(Q, e^{2x^2 Q}\right)$ for $Q = p(1-p)$ with $p$ uniform on (0,1), since $\mathbb{E}Q = 1/6$. This covariance is strictly positive for every $x > 0$. Scaling at negative infinity gives $U_0 \to 1/2$; completing the square and using $2A' = 1 + \lambda A$ gives $U_0 \to 0$ at positive infinity. Thus $U_0$ has the asserted analytic inverse.

Fix $u$ and bracket $\lambda_0(u)$ between finite values $\lambda_- < \lambda_+$. Continuity of (22) and strict decrease in $\lambda$ imply $\lambda_- < \lambda_u(d) < \lambda_+$ for small positive $d$. Letting the bracket shrink proves $\lambda_u(d) \to \lambda_0(u)$; the implicit function theorem then gives analytic continuation of the branch through $d = 0$. Hence $V(d) \to v_+(u)$.

At the other endpoint, for every finite $M$, positive variance implies $u(L, M) > L/[2(L+1)] = u$. It follows that $\lambda_u(d) \to +\infty$ as $d \uparrow L$. The uniform concentration estimate from Section 3.1 gives $J_0 J_2 / J_1^2 \to 1$ and $J_0^2 J_3 / J_1^3 \to 1$, so

$$\lim_{d \uparrow L} V(d) = \frac{L^2}{6(L+1)(L+2)} = v_-(u). \qquad (24)$$

Strict decrease in (20) implies $v_-(u) < v_+(u)$ and maps the entire interval $(0, L)$ bijectively onto $\left(v_-(u), v_+(u)\right)$. Equations (21) and the analytic inverse function theorem give the global analytic inverse. Amplitude and scale follow from $b_0 = e^a$ and $b_1 = e^a r / \left(d\sqrt{\alpha}\right)$, proving (9). Holding the first three coefficients fixes $u, b_0, b_1$ and leaves exactly this full analytic curve, proving the last assertion of Theorem 1.

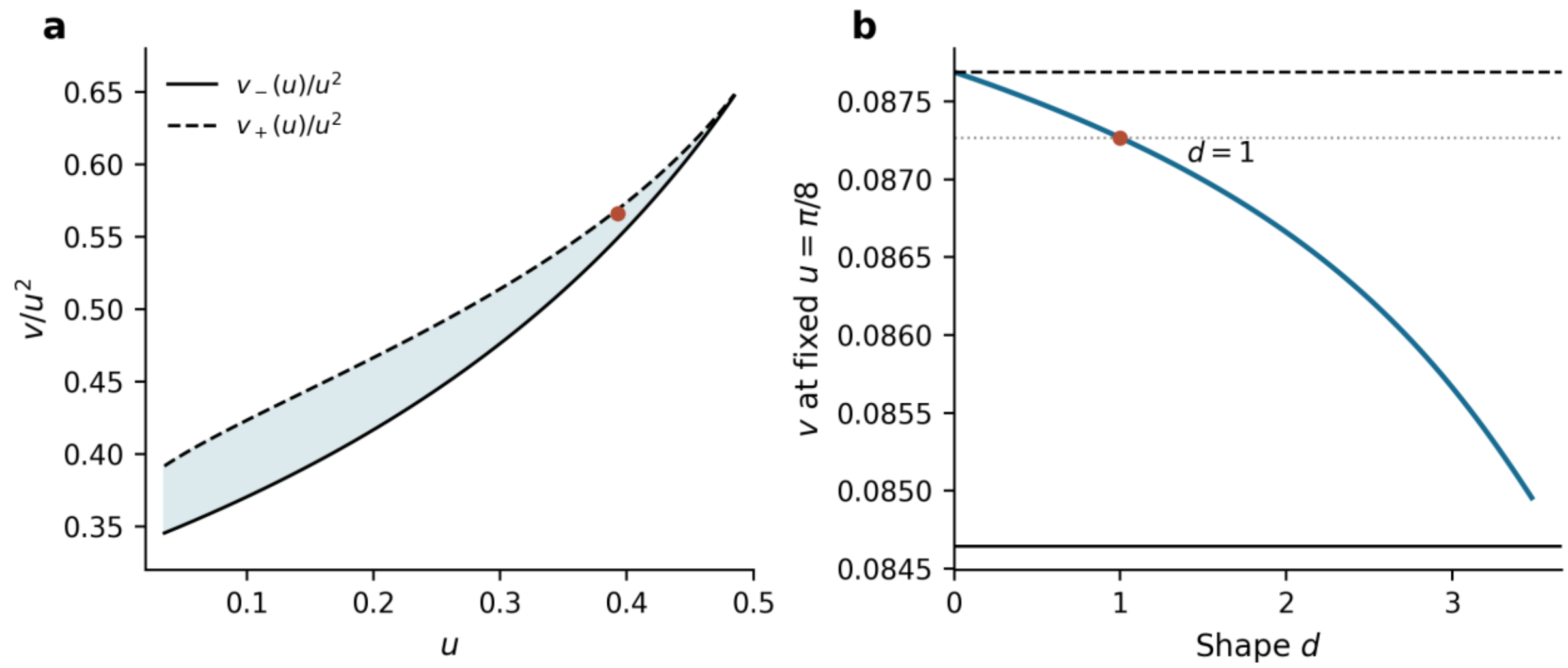

Figure 1. Quadratic shape geometry. Left: the realizability band in coordinates $(u, v/u^2)$, bounded by $v_-/u^2$ (solid) and $v_+/u^2$ (dashed). Right: the strictly decreasing fourth-coefficient ratio along the compensation curve with $u = \pi/8$. The marked point is $d = 1, \lambda = 0$. The curves are numerical illustrations; Theorem 1 establishes the full domain and monotonicity analytically.

## 4 Exact shape ranges for coefficient intervals

Theorem 1 concerns exact data. For the quadratic family, its geometry also gives a deterministic uncertainty statement that retains the dependence between $u$ and $v$ induced by (4).

**Theorem 2.** Let $\mathcal{B} = \prod_{i=0}^{3}[\ell_i, h_i]$ with $0 < \ell_i \leq h_i < \infty$, and define

$$b^- = (\ell_0, h_1, \ell_2, h_3), \qquad b^+ = (h_0, \ell_1, h_2, \ell_3). \qquad (25)$$

Every datum in $\mathcal{B}$ is realizable if and only if $b^-$ and $b^+$ are both realizable. Under this condition,

$$\{d(b): b \in \mathcal{B}\} = [d(b^-), d(b^+)]. \qquad (26)$$

The endpoints differ strictly whenever at least one coefficient interval has positive width.

**Proof.** Both boundary functions have elasticity greater than two:

$$\frac{uv_-{}'(u)}{v_-(u)} = 2 + \frac{u}{1-u} > 2, \qquad \frac{uv_+{}'(u)}{v_+(u)} > 2. \qquad (27)$$

For the second inequality, put $r_1 = A'/A$, $r_2 = A''/A'$, $r_3 = A'''/A''$. Equation (23) gives $c_0 = 2r_1 - r_2 > 0$. Equation (11), applied at shape one, gives $c_1 = 2r_2 - r_1 - r_3 > 0$. Taking logarithmic derivatives of (7) yields elasticity $2 + c_1/c_0 > 2$.

Let $x_i = \log b_i$. Then $\log u = x_0 - 2x_1 + x_2$ and $\log v = 2x_0 - 3x_1 + x_3$. If the two corners are realizable, $u(b^+) < 1/2$ ensures $u < 1/2$ throughout the box. Write $E_\pm = uv_\pm{}'/v_\pm$. The upper margin $\log v_+(u) - \log v$ has gradient $(E_+ - 2, 3 - 2E_+, E_+, -1)$, with signs $(+,-,+,-)$. The lower margin $\log v - \log v_-(u)$ has the opposite pattern. Therefore their minima on the box occur at $b^-$ and $b^+$ respectively, and both are positive. This proves full realizability before using the inverse map.

For its shape derivative, set $A_1 = 2R_1 - R_0 - R_2 > 0$ and $A_2 = 2R_2 - R_1 - R_3 > 0$. Then $u_\lambda = -uA_1$ and $v_\lambda = -v(2A_1 + A_2)$. With $J = \det DH < 0$, the chain rule gives

$$\nabla_x d = \frac{uv}{J}(-A_2,\ A_1 + 2A_2,\ -2A_1 - A_2,\ A_1). \qquad (28)$$

Its signs are $(+,-,+,-)$. Coordinatewise monotonicity gives the bounds (26), and continuity on the connected box gives every intermediate value. Strictness follows along any nondegenerate coordinate. Necessity of the corner test is immediate. ▫

The theorem tests whether **every** datum in the box has a model explanation. Failure of a corner test still allows some interior data to be realizable. Equation (26) is an exact marginal interval for $d$, not a rectangular description of the joint parameter set. Likewise, a deterministic coefficient box is not a confidence region unless a sampling model supplies its coverage.

For a concrete example, the quadratic model with normalized coefficients

$$(b_0, b_1, b_2, b_3) = (1, 1, \pi/8, \pi/36) \qquad (29)$$

has $d = 1$, $\lambda = 0$, $\alpha = 1/\pi$, and $\beta = 0$. Its first three coefficients admit all shapes $0 < d < \pi/(4-\pi)$. The fourth chooses $d = 1$. Figure 2 applies Theorem 2 to small independent relative intervals around (29). Scaling all four coefficients up or down together leaves both shape ratios unchanged; those two ordinary corners therefore cannot replace (25).

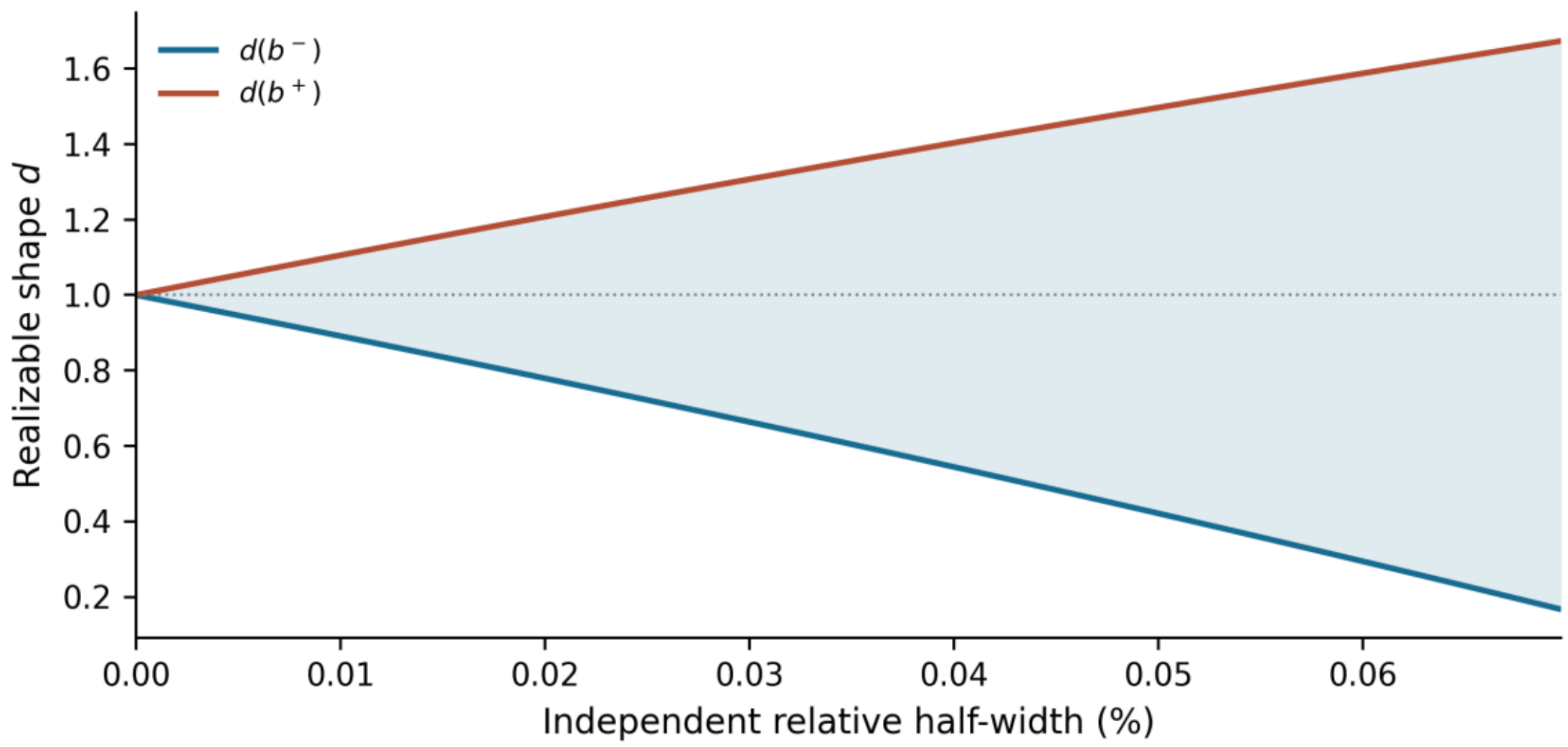

Figure 2. Minimum and maximum shape obtained by inverting the two alternating corners at each relative half-width. The shaded region is the interval characterized by Theorem 2. The plotted endpoint evaluations use numerical special functions and are illustrative rather than interval-certified values.

## 5 A global finite bound for every integer potential

**Theorem 3.** For every fixed integer $q \geq 2$, the map

$$\mathcal{B}_q : \Theta \to \mathbb{R}^{q(q+1)}, \qquad \theta \mapsto \left(b_0, \dots, b_{q(q+1)-1}\right) \qquad (30)$$

is injective and has differential rank four at every point. At each parameter point, four coefficients selected from this vector give local real analytic coordinates. The selected indices may depend on the point.

**Proof of injectivity.** Set $N = q(q+1) - 1$. Suppose two parameter sets have equal coefficients through $b_N$. Their amplitudes agree. Write $s_k = k!\, e^{-a} b_k > 0$ and

$$A_i(x) = \prod_{j=1}^{q-1} (x + d_i + j), \qquad q\alpha_i A_i(k) s_{k+q} = s_k + \beta_i s_{k+1} \quad (0 \leq k \leq N - q). \qquad (31)$$

This follows by integrating the derivative of $y^{d_i+k} e^{-\alpha_i y^q + \beta_i y}$, whose endpoint values vanish, and dividing the resulting moment recurrence by $(d_i)_{k+1}$.

Define $\delta = \beta_1 - \beta_2$, $P = \alpha_1 A_1 - \alpha_2 A_2$, and $Q = \beta_1 \alpha_2 A_2 - \beta_2 \alpha_1 A_1$. If $\delta = 0$, subtraction in (31) gives $P(k) = 0$ at at least $q$ distinct integers. Since $\deg P \le q - 1$, $P$ vanishes identically. Its leading two coefficients give $\alpha_1 = \alpha_2$ and $d_1 = d_2$.

If $\delta \neq 0$, subtraction and substitution give

$$qP(k)s_{k+q} = \delta s_{k+1}, \qquad P(k)s_k = Q(k)s_{k+1}. \qquad (32)$$

Positivity ensures that $P(k)$ and $Q(k)$ are nonzero at each used index. Multiplying the adjacent ratios in the second identity and comparing with the first shows that

$$R(x) = q\prod_{j=0}^{q-1} P\,(x+j) - \delta \prod_{j=1}^{q-1} Q\,(x+j) \qquad (33)$$

vanishes at $k = 0, \dots, N - 2q + 1$, a total of $q(q-1) + 1$ distinct integers. Yet its degree is at most $q(q-1)$ and it is nonzero. To establish the latter point, let $r = \deg P$ and $s = \deg Q$. An identity $R = 0$ would require $qr = (q-1)s$. With $0 \le r, s \le q - 1$, coprimality forces $r = s = 0$. This contradicts

$$Q + \beta_1 P = \delta \alpha_1 A_1, \qquad (34)$$

whose right side has degree $q - 1 \ge 1$. The contradiction excludes $\delta \neq 0$ and proves injectivity.

**Proof of full rank.** Let a tangent vector annihilate all coefficients in (30). Since $b_0 = e^{a}$, $\dot{a} = 0$, and all $\dot{s}_k = 0$. At the base point put

$$A(x) = \prod_{j=1}^{q-1}(x+d+j)\,, \quad P(x) = \dot{\alpha}A(x) + \alpha\dot{d}\ \partial_d A(x), \quad \delta = \dot{\beta}, \quad Q(x) = \delta\alpha A(x) - \beta P(x). \qquad (35)$$

Differentiating (31) gives precisely (32). If $\delta = 0$, then $P = 0$, whose leading coefficients give $\dot{\alpha} = \dot{d} = 0$. If $\delta \neq 0$, the same polynomial (33) has too many zeros; the nonzero argument now uses $Q + \beta P = \delta\alpha A$. Thus the tangent vector is zero. Analyticity was established in Section 2, and the analytic inverse function theorem applied to a nonzero four-row minor proves the local-coordinate statement. ▫

The proof actually identifies positive sequences satisfying the two displayed recurrences. It does not need an infinite moment-determinacy argument. Conversely, an arbitrary positive finite sequence satisfying a candidate recurrence need not equal the integrals of (1). Identification on the model image and realizability of arbitrary input remain separate questions.

## 6 Smaller coefficient sets and their exceptional cases

### 6.1 A finite-difference identification criterion

For data through $b_{2q}$ define, for $0 \le k \le q$,

$$f_k = \frac{(k+q)!}{k!} b_{k+q}, \qquad r_k = \frac{(k+1)b_{k+1}}{f_k}, \qquad t_k = \frac{b_k}{f_k}. \qquad (36)$$

The recurrence is equivalent to

$$A(k) - \beta r_k = t_k, \qquad A(x) = q\alpha \prod_{j=1}^{q-1}(x+d+j)\,. \qquad (37)$$

Treating $A$ as a polynomial of degree at most $q - 1$ gives $q + 1$ linear unknowns together with $\beta$.

**Proposition 4.** If $\Delta^q r_0 \neq 0$, the first $2q+1$ coefficients have at most one full parameter preimage. For realizable data,

$$\beta = -\frac{\Delta^q t_0}{\Delta^q r_0}, \qquad \alpha = \frac{[x^{q-1}]A}{q}, \qquad d = \frac{[x^{q-2}]A}{(q-1)[x^{q-1}]A} - \frac{q}{2}, \tag{38}$$

where $A$ is interpolated from $t_k + \beta r_k$ and $a = \log b_0$.

Indeed, a difference of two solutions to (37) satisfies $P(k) - c r_k = 0$, with $\deg P \leq q-1$. Taking a $q$th difference gives $c = 0$, then $P = 0$. The leading coefficients of $A$ give (38). An arbitrary input additionally requires positive recovered $d, \alpha$, the correct factorization of $A$, and agreement with the defining integral coefficients; the linear system alone certifies none of these existence conditions.

**Theorem 5.** For every fixed $q \geq 2$, every model with $\beta = 0$ satisfies $\Delta^q r_0 < 0$. Its first $2q+1$ coefficients therefore identify it against the whole tilted family, allowing competing $\beta \neq 0$. The condition $\Delta^q r_0 \neq 0$ holds on an open dense subset of $\Theta$.

**Proof.** At zero tilt, $I(s) = q^{-1}\alpha^{-s/q}\Gamma(s/q)$ and $r_k = (s)_{q-1} I(s)/I(s+q-1)$ with $s = d+k+1$. Set

$$\eta = 1/q, \quad x = (s+q-1)/q, \quad P(x) = \prod_{i=1}^{q-2}(x - i/q), \quad f(x) = \Gamma(x+\eta)/\Gamma(x). \tag{39}$$

Then $r_k = q^{q-1}\alpha^{(q-1)/q} P(x) f(x)$. The empty product for $q = 2$ is one. A beta integral, integrated by parts, gives

$$f(x) = C_\eta \int_0^\infty (1 - e^{-xt}) \frac{e^{-\eta t}}{(1-e^{-t})^{1+\eta}}\, dt = C_\eta \sum_{n=0}^\infty c_n \left(\frac{1}{z_n} - \frac{1}{x+z_n}\right), \tag{40}$$

where $C_\eta = \eta/\Gamma(1-\eta)$, $c_n = (1+\eta)_n/n!$, and $z_n = n + \eta$. In the integral, the change of variable $u = e^{-t}$ reduces integration by parts to the derivative of $\left(u/(1-u)\right)^\eta$; both boundary terms vanish. The positive series follows by expanding $(1-e^{-t})^{-1-\eta}$. Its terms, kept in the combined form shown, are $O(n^{\eta-2})$ locally uniformly for $x > 0$. This is the classical Bernstein structure of the gamma ratio; related broader families appear in [11].

Take the $q$th difference in $s$, corresponding to step $1/q$ in $x$. Polynomial division of $P(x)/(x+z)$ and the identity

$$\Delta_h^q \frac{1}{x+z} = \frac{(-1)^q q!\, h^q}{\prod_{j=0}^q (x+z+jh)} \tag{41}$$

give

$$\Delta^q r_k = -\frac{q!}{q}\alpha^{(q-1)/q} C_\eta \sum_{n=0}^\infty \frac{c_n \prod_{i=1}^{q-2}(z_n + i/q)}{\prod_{j=0}^q (x + z_n + j/q)} < 0. \tag{42}$$

The polynomial terms disappear, and the sign follows from $P(-z) = (-1)^{q-2}\prod_i (z + i/q)$. The resulting series has terms $O(n^{\eta-3})$ and converges absolutely. Finite differences can be exchanged with (40) because only finitely many function values occur. This proves the zero-tilt assertion. Finally, $\Delta^q r_0$ is real analytic on connected $\Theta$ and is nonzero on this slice, so its zero set has empty interior. Its nonzero set is open and dense, proving the generic statement. ▫

For $q = 3$, seven coefficients suffice whenever this third difference is nonzero. The theorem leaves open whether it can vanish at a real model with nonzero tilt. This exceptional-case question is unnecessary for the unconditional twelve-coefficient guarantee of Theorem 3.

### 6.2 A two-preimage bound

**Proposition 6.** For every fixed integer $q \geq 3$, any first $q + 3$ coefficients have at most two full parameter preimages in $\Theta$.

**Proof.** Normalize as in (4) and relabel the normalized coefficients as $b_k$, so $b_0 = b_1 = 1$ and $m_1 = d$. Let $u = b_2$, $v = b_3$, and $h = 1 - 2u$. The recurrence

$$q\alpha m_{k+q} = (d + k)m_k + \beta m_{k+1} \qquad (43)$$

and strict increase of moment ratios imply $h > 0$. Explicitly, under the probability measure tilted by $y$, the covariance of $Y$ and $Y^{q-1}$ is positive; hence $m_1 m_{q+1} > m_2 m_q$. Comparing (43) at $k = 0,1$ yields $(d + 1)m_1/m_2 > d/m_1$, equivalent to $u < 1/2$.

Put

$$\Delta_q(d) = (q + 1)(d + q)b_{q+1} - 2u(d + 1)b_q. \qquad (44)$$

The same first two recurrences yield

$$\alpha(d) = \frac{h}{q\,q!\,\prod_{j=2}^{q-1}(d + j)\ \Delta_q(d)}, \qquad \beta(d) = q\alpha(d)q!\prod_{j=1}^{q-1}(d + j)\,b_q - 1. \qquad (45)$$

Every genuine preimage has $\Delta_q(d) > 0$. Substitution into the $k = 2$ recurrence gives the quadratic or lower-degree polynomial

$$\begin{aligned} F_q(d) = {}& \frac{h(q + 1)(q + 2)}{2}(d + q)(d + q + 1)b_{q+2} \\ & - (u - 3v)(d + 2)\Delta_q(d) - 3vh(d + 1)(d + 2)b_q = 0. \end{aligned} \qquad (46)$$

It is not identically zero because its polynomial evaluation at $d = -2$ equals $h(q + 1)(q + 2)(q - 2)(q - 1)b_{q+2}/2 > 0$. No negative-shape density is used in this evaluation. There are at most two positive roots, and each determines $\alpha, \beta$ uniquely through (45). Undoing normalization recovers amplitude and scale uniquely. ▫

These polynomial roots are candidates. To verify actual realization for general $q$, one must also check the first $q - 1$ moments of the recovered density, then use (43) to propagate them. The proposition gives an upper bound of two, without asserting that two realized preimages exist at this coefficient count. When $q = 3$, it gives the six-coefficient bound in Table 1.

## 7 A certified failure of four-coefficient identification for the cubic potential

For $q = 3$, let $T$ have density proportional to $t^{d-1}e^{-t^3+\lambda t}$, let $r = \mathbb{E}T$, and normalize $Y = dT/r$. Then

$$a = 0, \quad \alpha = (r/d)^3, \quad \beta = \lambda r/d, \quad b_k = \frac{d^k \mathbb{E}T^k}{r^k k!\,(d)_k}, \qquad b_0 = b_1 = 1. \qquad (47)$$

Write $H_3(d, \lambda) = (b_2, b_3)$ and $w_3(d, \lambda) = b_4$ under this normalization.

**Theorem 7.** There are two legal parameter sets with shapes $d = 1$ and $d = 3$ whose first four coefficients agree exactly. More precisely, the equation

$$\mathcal{F}(\lambda_1, \lambda_3) = H_3(1, \lambda_1) - H_3(3, \lambda_3) = 0 \qquad (48)$$

has exactly one root in the maximum-norm box of radius $10^{-7}$ centered at

$$c = (-1.322527110062,\ 3.054536213083). \qquad (49)$$

Its distance from $c$ is less than $2^{-40}$. At this exact root, the fifth coefficients satisfy

$$\frac{10803884}{10^{12}} < w_3(1, \lambda_1) - w_3(3, \lambda_3) < \frac{10803885}{10^{12}}. \qquad (50)$$

**Proof.** The exact rational integration and interval arithmetic in Appendix A and the accompanying executable certificate establish, for

$$C = \begin{pmatrix} -752.30625069 & 1238.52173501 \\ -3473.96769536 & 6022.75907656 \end{pmatrix}, \qquad (51)$$

that $\det C \neq 0$, $\sup_X \parallel I - CD\mathcal{F} \parallel_\infty < 1/250$, and $\parallel C\mathcal{F}(c) \parallel_\infty < 6 \times 10^{-13}$. All displayed decimals in (49),(51) are exact rationals. Thus $x \mapsto x - C\mathcal{F}(x)$ maps the closed box $X$ strictly into itself and contracts. The contraction theorem gives its unique fixed point and therefore the unique zero in $X$. The a posteriori distance bound is less than $(6 \times 10^{-13})/(1 - 1/250) < 2^{-40}$. Evaluation on this smaller enclosure proves (50). Equation (47) gives the two full parameter sets, which differ because their shapes differ. ▫

The rerun certificate gives contraction bound $0.00395321869686$ and center-step bound $5.56064968456 \times 10^{-13}$, rounded upward for display. Table 2 gives approximate values for interpretation; the common exact data are defined by the root of (48), rather than by rounded decimal inputs.

| Shape d | α | β | $b_2$ | $b_3$ | $b_4$ |
|---|---|---|---|---|---|
| 1 | 0.05153096825 | −0.49214606557 | 0.40224172583 | 0.09125291848 | 0.013567817473 |
| 3 | 0.06848160832 | 1.24969278390 | 0.40224172583 | 0.09125291848 | 0.013557013588 |

Table 2. A pair of models sharing $b_0 = b_1 = 1$ and exactly the same $b_2, b_3$. Decimal values are approximations to the certified exact construction. The fifth coefficient separates this pair; it does not establish five-coefficient global identification.

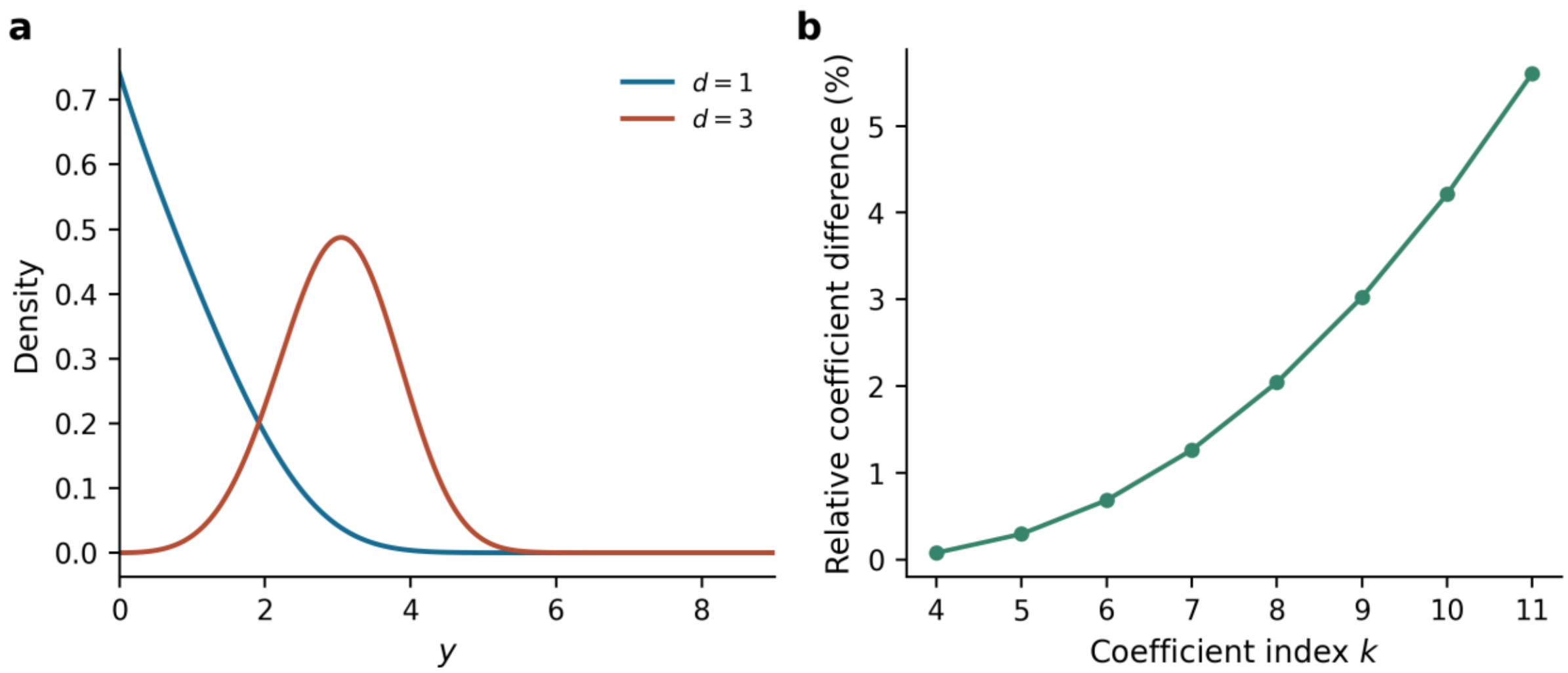

Figure 3. The certified cubic pair. Left: the two densities in the original normalized variable $Y$. Their ordinary means are $d = 1$ and $d = 3$, while $b_1 = 1$ for both because its denominator also contains $d$. Right: the relative difference $\left(b_k^{(1)} - b_k^{(3)}\right)/b_k^{(3)}$ for $k = 4, \ldots, 11$, evaluated numerically at the root. Only the existence, exact first-four equality, and fifth-coefficient gap are asserted by the rational certificate.

## 8 Identification and sensitivity

Full rank gives a local analytic inverse through suitable coefficient coordinates, and hence local Lipschitz dependence after restricting to a sufficiently small neighborhood. It does not supply a stability constant over the noncompact parameter domain. The following elementary boundary result makes this distinction precise.

**Proposition 8.** Fix $q \geq 2$ and any finite coefficient count $K$ for which the coefficient map is globally injective. Its inverse on the image is not uniformly continuous when parameter and data spaces have their usual Euclidean metrics.

**Proof.** Fix two shapes $d_1 \neq d_2$, take $a = 0$, $\beta = -1$, and let the common positive $\alpha \downarrow 0$. Dominated convergence in (1), using $y^{d_i+k-1}e^{-y}$, gives $m_k \to (d_i)_k$. Thus both coefficient vectors converge to $(1,1,1/2!,\dots,1/(K-1)!)$, while the parameter vectors remain at least $|d_1 - d_2|$ apart. This contradicts uniform continuity of the inverse. ▫

At the excluded boundary $\alpha = 0$, every shape produces the same full normalized coefficient sequence at fixed $a$ and negative $\beta$. This explains why exact identification in the open domain can coexist with severe loss of shape sensitivity near its boundary. The coefficient-box theorem supplies a sharp deterministic treatment in the quadratic family; Proposition 4 exposes a different source of instability through a small finite-difference denominator in higher potentials.

The results also delimit the role of a continuous dimension interpretation. With $D = 2d$, the quadratic theorem says that a dimension-like change can be compensated exactly by scale and tilt at three consecutive coefficients, and that the fourth breaks this compensation everywhere in the stated model. The cubic construction shows that this threshold depends on the potential. Extending these statements to unknown potential degree, additional potential terms, or a physical measurement process requires further work.

## 9 Conclusions

The Pochhammer normalization makes continuous shape identification a distinct finite-data problem. For a quadratic potential, four consecutive transform coefficients give a complete global inverse and an exact realizability region; three leave an analytic compensation curve. The shape geometry also determines exact coefficient-box extrema from two alternating corners. For any fixed integer potential, recurrence comparison gives an unconditional finite identification bound and an everywhere full-rank map, while a finite-difference criterion supplies smaller generic bounds. A rationally certified cubic pair establishes a genuine failure of the quadratic four-coefficient threshold. These statements provide a self-contained basis for studying how model structure controls the amount of response information needed to separate shape, scale, and tilt.

## Data and code availability

The accompanying reproducibility package contains the rational cubic moment evaluator, the exact collision certificate and its rerun output, symbolic checks of the displayed elimination identities, and scripts and numerical data for all figures. The certificate uses Python integer and rational arithmetic. Figure calculations use NumPy, SciPy, mpmath, and Matplotlib and are marked as illustrations. All mathematical inputs needed for the stated theorems are specified in the article and Appendix A.

## Appendix A Rational certification of the cubic collision

For integer $n \geq 0$, write $I_n(\lambda) = \int_0^\infty t^n\, e^{-t^3+\lambda t}\, dt$. The certificate uses $n \leq 8$ and rational $|\lambda| < 4$. Differentiation gives $I_n{}' = I_{n+1} > 0$, so an interval in $\lambda$ is enclosed by evaluations at its endpoints. Derivatives of the shape ratios are computed from

$$u_\lambda = u(R_0 + R_2 - 2R_1), \qquad v_\lambda = v(2R_0 + R_3 - 3R_1), \quad R_j = I_{d+j}/I_{d+j-1}. \qquad (A1)$$

Only gamma values at integer multiples of $1/3$ are needed. For $m = 0,1$, use $\Gamma\big((m+1)/3\big) = 3\int_0^\infty t^m\, e^{-t^3}\, dt$. Let $T \geq 1$ be an integer and $N$ even. The alternating Taylor inequalities for $e^{-x}$ at $x \geq 0$ give rational bounds on the integral over $[0, T]$ with upper sum

$$G_N(m,T) = 3\sum_{j=0}^{N} \frac{(-1)^j T^{3j+m+1}}{j!\,(3j+m+1)}. \qquad (A2)$$

The lower sum is $G_{N+1}$. For $t = T + y$, convexity gives $t^3 \geq T^3 + 3T^2 y$, so the omitted positive tails are bounded above by

$$B_0 = \frac{2^{-T^3}}{T^2}, \qquad B_1 = 3\left(\frac{T}{3T^2} + \frac{1}{9T^4}\right)2^{-T^3}. \qquad (A3)$$

Thus the gamma values lie in $[G_{N+1}, G_N + B_m]$. The value $\Gamma(1) = 1$ and the recurrence $\Gamma\big((r+3)/3\big) = (r/3)\Gamma(r/3)$ supply the rest.

Expand the tilt to degree $L$:

$$I_n(\lambda) = \frac{1}{3}\sum_{j=0}^{L} \frac{\lambda^j}{j!}\Gamma\big((n+j+1)/3\big) + E_{n,L}. \qquad (A4)$$

Set $B = |\lambda|$ and choose an integer $R \geq 0$ with $R^2 \geq 2B$. For all $t \geq 0$, $Bt \leq t^3/2 + BR$. The exponential Taylor remainder, $e < 4$, and the gamma upper enclosure give

$$\left|E_{n,L}\right| \leq 2^{\lceil 2BR \rceil} \frac{B^{L+1}}{(L+1)!} \frac{2^{\lceil (n+L+2)/3 \rceil}}{3}\, \overline{\Gamma}\big((n+L+2)/3\big). \qquad (A5)$$

This bound applies to either sign of $\lambda$. Gamma base integrals are summed as exact fractions. Subsequent arithmetic rounds outward to a dyadic grid, using floor for lower endpoints and ceiling for upper endpoints; multiplication checks all endpoint products, and division first excludes zero in the denominator. Precision and truncation are increased until each required moment interval is positive and has width at most $2^{-80}$. The actual output records all truncation degrees, gamma bounds, moment intervals, and matrix inequalities, allowing the finite computation to be inspected separately from its implementation.

At the rational center and at the coarse and refined endpoints, the program performs ten distinct point evaluations. It then bounds the Jacobian on the whole coarse box by monotonicity of the raw $I_n$, establishes the contraction inequalities of Theorem 7, and encloses the fifth-coefficient difference on the refined box. The final rational comparisons include the two strict bounds in (50). This is a computer-assisted existence proof for one pair, with uniqueness restricted to the stated two-tilt box. It does not count all parameter preimages of the common data.